\documentclass[12pt]{amsart}

\usepackage{graphicx}
\usepackage{amsthm}
\usepackage{amsmath, amssymb}
\swapnumbers

\theoremstyle{plain}
\newtheorem{Thm}{Theorem}

\newtheorem{Coro}[Thm]{Corollary}
\newtheorem{Lem}[Thm]{Lemma}

\newtheorem{Prop}[Thm]{Proposition}

\theoremstyle{definition}

\begin{document}

\title{Mapping class groups of Heegaard splittings}

\author{Jesse Johnson and Hyam Rubinstein}
\address{\hskip-\parindent
        Department of Mathematics \\
        Oklahoma State University \\
        Stillwater, OK 74078 \\
        USA}
\email{jjohnson@math.okstate.edu}
\address{\hskip-\parindent
        Dept. of Mathematics and Statistics \\
        The University of Melbourne \\
        Parkville, VIC, 3053 \\
        AUSTRALIA}

\email{H.Rubinstein@ms.unimelb.edu.au}

\thanks{Research supported by NSF MSPRF grant 0602368 and ARC grant DP020849}

\begin{abstract}
The mapping class group of a Heegaard splitting is the group of automorphisms of the ambient 3-manifold that take the surface onto itself, modulo isotopies that keep the surface on itself. We characterize the mapping classes that restrict to periodic and reducible automorphisms of the surface.
\end{abstract}

\maketitle

Given a 3-manifold $M$, define $Aut(M)$ to be the group of automorphisms (self homeomorphisms) of $M$.  The set of connected components of $Aut(M)$ forms a group $Mod(M)$ called the \textit{mapping class group} of $M$.  Given a Heegaard splitting $(\Sigma, H^-, H^+)$ of $M$, we will define $Aut(M, \Sigma)$ to be the subset of $Aut(M)$ consisting of all maps that send $\Sigma$ onto itself.  The set of connected components of $Aut(M, \Sigma)$ forms a group, the \textit{mapping class group of $\Sigma$}, which we will denote $Mod(M, \Sigma)$.
Note that both the areas of automorphisms of $3$-manifolds (and knots and links in $S^3$) and Heegaard splittings of $3$-manifolds (tunnel systems for knots and links) are highly developed. 

There has been some work on the connection between these areas, as a tool to obtain bounds for the number of Heegaard splittings or the size of the automorphism group of a $3$-manifold, in terms of each other. However the natural connection between these areas, namely the study of the mapping class group $Aut(M, \Sigma)$ of a Heegaard splitting $\Sigma$, has not been systematically investigated. One key to understanding $Aut(M, \Sigma)$ is the observation that an element of the kernel of the map from $Aut(M, \Sigma)$ to $Aut(M)$ corresponds to a non-trivial loop of embeddings of the Heegaard surface. In fact, for a hyperbolic 3-manifold the space of embedded surfaces isotopic to $\Sigma$ is a classifying space for the kernel~\cite{JMc}. Just as strongly irreducible Heegaard surfaces define minimal surfaces of index one, such loops of embeddings should determine minimal surfaces of index two. See Bachman~\cite{bach} and Hass, Thompson, Thurston~\cite{htt} for further developments of similar ideas. 

An element of $Mod(M, \Sigma)$ will be called \textit{periodic, reducible or pseudo-Anosov} if its restriction to $\Sigma$ is periodic, reducible or pseudo Anosov, respectively.  We will present below characterizations of periodic and reducible elements of $Mod(M, \Sigma)$, as well as showing that stabilized Heegaard splittings admit a large number of pseudo-Anosov automorphisms. 

One can construct a Heegaard splitting with periodic elements in its mapping class group as follows:  Let $N$ be a 3-manifold and $M$ a finite, regular branched or unbranched cover of $N$. Let $(R, G^-, G^+)$ be a Heegaard splitting of $N$.  If the covering of $N$ is branched, assume the branch set $B$ is in bridge position with respect to $\Sigma$, i.e. that every component of $B \setminus \Sigma$ is an arc or tree that is parallel into $\Sigma$.  This condition on the branch set implies that the Heegaard splitting $(R, G^-, G^+)$ lifts to a Heegaard splitting $(\Sigma, H^-, H^+)$ of $M$.  

The deck transformations of the covering define a finite group of automorphisms of $M$.  These automophisms take $\Sigma$ onto itself so they define a finite subgroup of $Mod(M, \Sigma)$.  We prove the following in Section~\ref{extendsect}:

\begin{Thm}
\label{mainthm1}
Let $G \subset Mod(M, \Sigma)$ be a finite subgroup of the mapping class group of a Heegaard splitting $(\Sigma, H^-, H_+)$.  Then $M$ is a (possibly branched) finite cover of a manifold $N$ such that $(\Sigma, H^-, H^+)$ is the lift of a Heegaard splitting $(R, G^-, G^+)$ for $N$, the branch set of the covering is in bridge position with respect to $R$ and the group $G$ is induced by the deck transformations of the covering.
\end{Thm}

Each connected component of $Aut(M, \Sigma)$ is a subset of a connected component of $Aut(M)$, so the inclusion map in $Aut(M)$ determines a canonical homomorphism $i : Mod(M, \Sigma) \rightarrow Mod(M)$.  In general, this map need be neither injective nor surjective. Elements in the kernel of $i$ correspond to isotopy trivial automorphisms of $M$. In a hyperbolic 3-manifold, finite order automorphisms cannot be isotopy trivial, so Theorem~\ref{mainthm1} implies the following:

\begin{Coro}
\label{mainthm2}
Assume $M$ is a closed hyperbolic 3-manifold. Then the kernel of the canonical homomorphism $i : Mod(M,\Sigma) \rightarrow Mod(M)$ contains no finite order elements and $i$ is one-to-one if and only if $Mod(M, \Sigma)$ is finite.
\end{Coro}

Reducible elements of $Mod(M, \Sigma)$ defy any single classification such as that of periodic elements. On the one hand reducible automorphisms can have ``local'' behaviour defined by an isotopy within a proper submanifold of $M$.  In Section~\ref{wrsect}, we show that such an automorphism exists whenever there are disks $H^-$, $H^+$ that intersect in exactly two points.  This intersection need not be minimal, so in particular every weakly reducible Heegaard splitting contains such an element.  

The second type of reducible automorphism, described in Section~\ref{opensect}, comes from an isotopy that sweeps through the entire manifold and is defined by an open book decomposition. A third type is induced by a one-sided Heegaard splitting of $M$ or a structure we call a branched Heegaard splitting. We describe reducible automorphisms of such Heegaard splittings in Section~\ref{onesidesect}, then in Section~\ref{refsect}, we show that under reasonable conditions, every reducible automorphism is either local or induced by an open book decomposition, a one-sided Heegaard surface or a branched Heegaard surface. 

\begin{Thm}
\label{redthm}
Assume $M$ is hyperbolic and $\Sigma$ is strongly irreducible.  If $\phi \in Mod(M, \Sigma)$ is reducible and $L \subset \Sigma$ is the fixed set of $\phi$ then either $\phi$ is induced by an open book decomposition, a one-sided Heegaard splitting or a branched Heegaard splitting that defines $\Sigma$ or $\phi$ is periodic outside a submanifold $N \subset M$ whose boundary is a union of surfaces that result from maximally compressing subsurfaces of $\Sigma$ in the complement of $L$.
\end{Thm}

Finally, in Section~\ref{stabsect}, we consider the case of pseudo Anosov maps.  In this case, we are only able to show that every stabilized Heegaard splitting has pseudo Anosov elements in its kernel.  The problem of characterizing in general Heegaard splittings with pseudo Anosov automorphisms remains open.

\section{The Mapping Class Group}
\label{mcgsect}

A Heegaard splitting for a 3-manifold $M$ is a triple $(\Sigma, H^-, H^+)$ where $H^-, H^+ \subset M$ are handlebodies (manifolds homeomorphic to regular neighborhoods of graphs) and $\Sigma$ is a compact, connected, closed and orientable surface embedded in $M$ such that $H^- \cup H^+ = M$ and $\partial H^- = \Sigma = \partial H^+$.  

As noted above, $Mod(M, \Sigma)$ is the group of equivalence classes of automorphisms of $M$ that take $\Sigma$ onto itself.  Two automorphisms are equivalent if there is an isotopy from one to the other by automorphisms that take $\Sigma$ onto itself.  In this paper, we will restrict our attention to elements of $Mod(M, \Sigma)$ that are orientation preserving on $M$ and orientation preserving on $\Sigma$.  Such automorphisms take the handlebody $H^-$ onto itself and $H^+$ onto itself.  If automorphisms that swap the handlebodies or reverse the orientation on $M$ exist, the orientation preserving automorphisms that preserve the handlebodies form an index two or four subgroup.  

Throughout the paper, $i$ will refer to the homomorphism $Mod(M, \Sigma) \to Mod(M)$ induced by ``forgetting'' $\Sigma$ and considering each element of $Mod(M, \Sigma)$ as an automorphism of $M$.  The two immediately obvious questions to ask are when $i$ will be onto and when $i$ will be one-to-one.  The first of these questions relates to two classical questions from the field of Heegaard splittings.

Two Heegaard surfaces, $\Sigma$ and $R$ of a manifold $M$ are called \textit{homeomorphic} if there is an automorphism of $M$ taking $\Sigma$ onto $R$.  The surfaces are \textit{isotopic} if there is an automorphism that is isotopic to the identity and sends $\Sigma$ onto $R$.

For any Heegaard surface $R$ homeomorphic to $\Sigma$, there is an automorphism $\phi$ of $M$ such that $\phi(\Sigma) = R$.  If the homomorphism $i$ is onto then there is an automorphism $\psi$ that is isotopic to $\phi$ and sends $\Sigma$ onto itself.  The map $\phi \circ \psi^{-1}$ is isotopic to the identity on $M$ and sends $\Sigma$ onto $R$, so $\Sigma$ is in fact isotopic to $R$.  Combining this argument with its converse, we get the following:

\begin{Prop}
\label{ontoprop}
The homomorphism $i$ is onto if and only if every Heegaard splitting that is homeomorphic to $\Sigma$ is isotopic to $\Sigma$.
\end{Prop}

Li~\cite{li:wald} showed that an atoroidal manifold has only finitely many isotopy classes of Heegaard splittings for each genus. This implies that for atoroidal manifolds, the image of $i$ has finite index in $Mod(M)$.  

For toroidal manifolds, Morimoto and Sakuma~\cite{ms:tunnels} have found a number of tunnel number one knot complements and Bachman and Derby-Talbot~\cite{bdt:seif} have found Seifert fibered spaces with an infinite number of homeomorphic, non-isotopic Heegaard splittings.  For all these Heegaard splittings, the image of $i$ has infinite index.  The question of the kernel of $i$ appears to have no classical analogue.

If the automorphism of $\Sigma$ is isotopic to the identity then it is isotopic to the identity on each handlebody, so $\phi$ is trivial in $Mod(M, \Sigma)$.  Thus the restriction map from $Mod(M, \Sigma)$ to $Mod(\Sigma)$ is an injection.  An element of $Mod(M, \Sigma)$ will have finite order if and only if its restriction to $\Sigma$ has finite order.

Mapping class groups of Heegaard splittings have been studied in a number of different contexts.  As part of a study of tunnel number one knots in $S^3$, Goeritz~\cite{goer:s3} showed that the mapping class group of a genus two Heegaard splitting of $S^3$ is finitely generated. Scharlemann~\cite{schar:mcg} recently published a new proof of this result, then shortly afterwards, the group was shown to be finitely presented independently by Akbas~\cite{akbas} and Cho~\cite{cho}. Scharlemann's proof was also a main tool in Cho and McCullough's work on the tree of unknotting tunnels~\cite{mccho}. A finite generating set for the mapping class group of the standard genus-three Heegaard splitting of the 3-torus has been found by the first author~\cite{me:t3}.

Futer~\cite{fu:invol} showed that the only tunnel number one knots with an automorphism of the complement that reverses the orientation of the knot and preserves a genus two Heegaard splitting are two bridge knots.  The genus two Heegaard splittings are those coming from level tunnels joining the two maxima or the two minima.

In order to study automorphisms of 3-manifolds, Birman and the second author have classified the mapping class groups of a number of one sided Heegaard splittings~\cite{rb:homeo}.  Zimmermann~\cite{zim:Hur} used automorphisms of Heegaard splittings to find hyperbolic 3-manifolds with large mapping class groups and Oertel~\cite{oer:hbdy} has given a nice characterization of automorphisms of handlebodies.  Oertel and Carvalho~\cite{oc:auts} used this to characterize automorphisms of reducible manifolds.  

Casson and Long~\cite{cl:alg} found an algorithm to determine if a pseudo Anosov map on the boundary of a handlebody extends to the entire handlebody.  They did this not in the context of Heegaard splittings, but to determine when a knot is not homotopically ribbon.  Long later found the first (and so far only) example of an irreducible Heegaard splitting with a pseudo Anosov automorphism~\cite{long:twoh}.

The first author and Darryl McCullough showed that the homotopy type of the space of embeddings of (unmarked) surfaces isotopic to $\Sigma$ is determined by a certain long exact sequence involving $Mod(M, \Sigma)$ and the homotopy groups of $Diff(\Sigma)$ and $Diff(M)$. For $M$ a hyperbolic 3-manifold, the space of embeddings turns out to be a classifying space for the kernel of the homomorphism $i : Mod(M, \Sigma) \rightarrow Mod(M)$, i.e. a space with this group as its fundamental group and a contractible universal cover.

A Heegaard splitting $(\Sigma, H^-, H^+)$ determines two subsets $\mathcal{H}^-$, $\mathcal{H}^+$ in the curve complex $\mathcal{C}(\Sigma)$ consisting of the loops bounding disks in $H^-$, and $H^+$ respectively. The edge-path distance between these two sets is the \textit{distance} of the Heegaard splitting, as defined by Hempel~\cite{hempel}. For genus three and greater, the group $Mod(M, \Sigma)$ is isomorphic to the group of isometries of $\mathcal{C}(\Sigma)$ that take each set $\mathcal{H}^\pm$ onto itself.  Namazi~\cite{nam:mcg} used this fact to show that for sufficiently high distance (depending on the genus) Heegaard splittings, the mapping class group of the Heegaard splitting is finite. The first author later showed using purely topological methods that this is true for any Heegaard splitting with distance four or greater in any genus~\cite{me:medium}. 

More recently, Campisi and Rathbun have shown that the disk sets $\mathcal{H}^\pm$ are determined by their large scale geometry~\cite{rathcam}. They conclude that $Mod(M, \Sigma)$ is isomorphic to the group of quasi-isometries of $\mathcal{C}(\Sigma)$ (modulo bounded-distance maps) that keep each set $\mathcal{H}^\pm$ within a bounded neighborhood of itself. Thus topological results about the mapping class group can be translated into results about quasi-isometry groups in the complex of curves.

One can also deduce a finite mapping class group from Lustig and Moriah's double rectangle condition~\cite{lm:dblrect}.  If a Heegaard splitting has one or more diagrams satisfying the double rectangle condition then these are permuted by the mapping class group, and Lustig and Moriah have shown that every Heegaard splitting has at most finitely many such diagrams.

\section{Periodic automorphisms}
\label{extendsect}

To prove Theorem~\ref{mainthm1} for a finite subgroup $G$ of $Mod(M, \Sigma)$, we would like to find an isomorphic subgroup $G'$ of $Aut(M, \Sigma)$ that maps onto $G$ by the inclusion map.  The Nielsen Realization Theorem ~\cite{ker:niels} implies that there is a subgroup $G''$ of $Aut(\Sigma)$ with this property.  The Equivariant Disk Lemma in~\cite{meeksyau} states that there is a system of disks in each handlebody that behaves nicely under $G''$, and we will use this system of disks to extend $G''$ into the handlebodies.

We paraphrase the Equivariant Disk Lemma as follows:

\begin{Thm}[Equivariant Disk Theorem - Theorem 7 in~\cite{meeksyau}]
\label{permutediskslem}
Let $G'' \subset Isom(\Sigma)$ be a finite group of isometries of $\Sigma = \partial H$ such that each isometry extends to an automorphism of $H$.  Then there is a collection $L$ of essential, simple closed geodesics in $\Sigma$, bounding pairwise disjoint, properly embedded disks such that $G''$ permutes the loops in $L$ and the complement in $\Sigma$ of $L$ is planar.
\end{Thm}

\begin{Coro}
\label{extendhlem}
Let $G'' \subset Isom(\Sigma)$ be a finite group of isometries of $\Sigma = \partial H$ such that each isometry extends to an automorphism of $H$.  Then there is an isomorphic subgroup $G' \subset Aut(H)$ such that each element of $G'$ restricts to an element of $G''$ and the set of fixed points form boundary parallel arcs and trees in $H$.
\end{Coro}

\begin{proof}
We will construct the group $G'$ in the piecewise linear category. Let $L$ be a collection of simple closed curves, guaranteed by the Equivariant Disk Theorem, that are permuted by $G''$ and bound a collection $\mathbf{D}$ of disks in $H$.  Choose a triangulation of $\Sigma$ such that $L$ is contained in the 1-skeleton, each isometry of $G''$ is simplicial and each fixed point is a vertex of the triangulation.

Extend this triangulation to each disk by taking the cone over the triangulation of its boundary. This triangulation consists of a vertex in the center of the disk and edges from this vertex to each vertex in the boundary.  Then for each $g \in G''$, the cone structure induces a unique simplicial map that fixes the interior vertex and agrees with $g$ on the boundary of the disk.

The complement of the disks is a collection of balls whose boundaries are now triangulated.  Extend this triangulation to the balls by taking a cone over the boundary of each ball. We can again extend $g$ to $H$ by the unique simplicial map that fixes the interior vertex and agrees with $g$ on the boundaries of the balls.  

For each element $g$ of $G''$, there is a unique simplicial map of $H$ that agrees with $g$ on the boundary.  The set $G' \subset Aut(H)$ consisting of these maps is a subgroup such that each element of $G'$ restricts to an element of $G''$.  The fixed point set in $H$ consists of the edges that have one vertex at a fixed point in $\partial H$ and the other vertex in the interior of $H$, or vertices at the center of a disk and a ball in the interior.  Any properly embedded arcs and graphs consisting of these edges are boundary parallel.
\end{proof}

\begin{proof}[Proof of Theorem~\ref{mainthm1}]
Let $G \subset Mod(M, \Sigma)$ be a finite subgroup of the mapping class group of a Heegaard splitting $(\Sigma, H^-, H^+)$.  The restriction of $G$ to $\Sigma$ is a subgroup of the mapping class group $Mod(\Sigma)$.  By the Nielsen realization Theorem~\cite{ker:niels}, there is a hyperbolic structure on $\Sigma$ and a subgroup $G''$ of $Isom(\Sigma)$ such that the restriction map from $G$ to $G''$ is an isomorphism.

By Corollary~\ref{extendhlem}, the group $G''$ extends to a finite subgroup of $Aut(H^-)$ and to a finite subgroup of $Aut(H^+)$.  Thus $G''$ extends to a subgroup $G'$ of $Aut(M)$.  

Let $N$ be the quotient of $M$ by orbits of $G'$.  Because $G''$ is piecewise-linear and finite, the inclusion map $M \rightarrow N$ is a (possibly branched) finite cover and $G''$ is the group of deck transformations for this cover.  Because the set of fixed points is boundary parallel in each handlebody, the preimage in $N$ of each handlebody is a handlebody, the Heegaard surface descends to a Heegaard surface $R$ and the branch set is in bridge position with respect to $R$.
\end{proof}

Theorem~\ref{mainthm1} allows us to pose the question of periodic elements of the kernel in terms of deck transformations of finite covers.  Or, equivalently, we can think of the deck transformations as a group action on $M$ by the finite subgroup of $Mod(M,\Sigma)$.  This allows us to employ a recent Theorem of Dinkelbach and Leeb.  The following is Theorem 1.2 from~\cite{dl:ricci}:

\begin{Thm}[Dinkelbach and Leeb~\cite{dl:ricci}]
\label{trivisolem}
Any smooth action by a finite group on a closed hyperbolic 3-manifold is smoothly conjugate to an isometric action.
\end{Thm}

\begin{proof}[Proof of Corollary~\ref{mainthm2}]
Let $(\Sigma, H^-, H^+)$ be a Heegaard splitting of a hyperbolic 3-manifold $M$ and let $\phi$ be a periodic element of $Mod(M, \Sigma)$. Let $G$ be the (finite) subgroup of $Mod(M, \Sigma)$ generated by $\phi$.  By Theorem~\ref{mainthm1}, $M$ is a (possibly branched) finite cover of a manifold $N$ such that $\phi$ is induced by a deck transformation of the covering.  The deck transformations induce a smooth action of $M$ by $G$, so Theorem~\ref{trivisolem} implies that the action by $G$ is conjugate to an action by isometries.  No isometry of a hyperbolic 3-manifold is isotopy trivial, so $\phi$ is conjugate to an isotopy non-trivial automorphism of $M$.  This implies that $\phi$ is isotopy non-trivial, and therefore $\phi$ is not in the kernel of the homomorphism $i : Mod(M,\Sigma) \rightarrow Mod(M)$.
\end{proof}

\section{Weakly Reducible Splittings}
\label{wrsect}

A Heegaard splitting $(\Sigma, H^-, H^+)$ is \textit{weakly reducible} if there are properly embedded, essential disks $D_1$, $D_2$ in $H^-$, $H^+$, respectively such that the boundaries of $D_1$ and $D_2$ are disjoint in $\Sigma$.  In this section we show that for every weakly reducible Heegaard splitting, the kernel of the induced homomorphism $i$ contains infinite order, reducible elements.  However, we begin by showing that a slightly weaker condition is enough to imply this.

\begin{Lem}
\label{twopointslem}
Let $(\Sigma, H^-, H^+)$ be a Heegaard splitting of genus at least 3.  If there are essential, properly embedded disks $D_1 \subset H^-$ and $D_2 \subset H^+$ whose boundaries intersect transversely in precisely two points then the kernel of $i$ contains an infinite order reducible element or the boundary loops of $D_1$ and $D_2$ are isotopic.
\end{Lem}

\begin{proof}
Let $(\Sigma, H^-, H^+)$ be a Heegaard splitting of a manifold $M$ of genus at least 3 and let $D_1 \subset H^-$, $D_2 \subset H^+$ be properly embedded, essential disks whose boundaries are not isotopic in $\Sigma$.   Assume $\partial D_1 \cap \partial D_2$ consists of two points.  If we orient $\partial D_1$ and $\partial D_2$ then at each intersection the orientations of the loops will induce an orientation of $\Sigma$.  We have two cases to consider: when the induced orientations agree and when they disagree.

Let $N$ be the closure of a regular neighborhood of $D_1 \cup D_2$.  Because there are two points of intersection, $N$ is a solid torus.  The portion of $\Sigma$ inside $N$ is shown, cut along a meridian, in Figure~\ref{localfig}. We can recover $\Sigma \cap N$ by gluing the top to the bottom, either directly or with a half twist.  If the orientations at the two intersections disagree then we glue directly and the surface $\Sigma \cap N$ will be a four punctured sphere whose boundary consists of four simple closed curves in $\partial N$. Two loops with this intersection pattern are shown on the left in Figure~\ref{disksfig} and the reader will note that a regular neighborhood is a four-times punctured sphere.  

If the orientations agree then we glue with a twist and $\Sigma \cap N$ will be a twice punctured torus, whose boundary consists of two simple closed curves in $\Sigma$. Two loops with this intersection pattern are shown on the right side of Figure~\ref{disksfig}, and the reader will see that they sit in a torus.  The loops $\Sigma \cap \partial N$ are four parallel longitudes of $N$ in the first case, and two loops that each wrap twice around a longitude of $N$ in the second case.
\begin{figure}[htb]
 \begin{center}
 \includegraphics[width=1.5in]{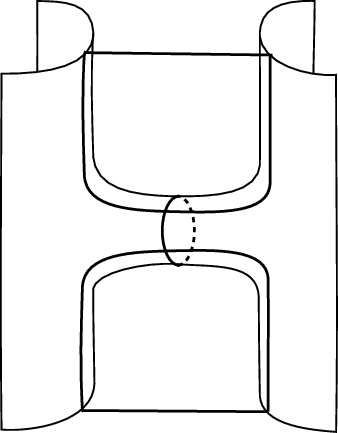}
\caption{A neighborhood of the two disks is a solid torus.}
  \label{localfig}
  \end{center}
\end{figure}
\begin{figure}[htb]
  \begin{center}
  \includegraphics[width=2.5in]{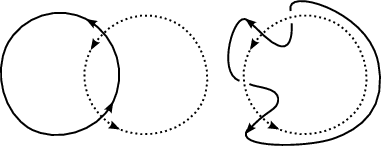}
  \caption{If disks intersect in two points, the induced orientations at the intersections may agree or disagree.}
 \label{disksfig}
 \end{center}
\end{figure}

In both cases, ``spinning'' $N$ along its longitude induces an automorphism of the solid torus $N$ that fixes the boundary of $N$ and induces an automorphism of $\Sigma \cap N$ consisting of Dehn twists along loops parallel to $\Sigma \cap \partial N$.  If the induced automorphism is trivial then each loop of $\Sigma \cap \partial N$ either bounds a disk in $\Sigma$ or cobounds an annulus with a second loop in $\Sigma \cap \partial N$.

In the case when the orientations disagree, there are four loops.  If two of the loops in $\Sigma \cap \partial N$ are trivial in $\Sigma$ then either one of the loops $\partial D_1$, $\partial D_2$ is trivial in $\Sigma$ or the two loops are parallel.  Because we assumed both loops are essential, we conclude that in this case, the two loops are parallel.

If each loop in $\Sigma \cap \partial N$ is parallel to a second loop in $\Sigma \cap \partial N$ then $\Sigma$ is the union of the four punctured sphere $\Sigma \cap N$ and two annuli.  In this case $\Sigma$ has genus two.  We assumed that $\Sigma$ has genus at least 3 so in the case when the orientations at the intersections agree, we conclude that the induced automorphism is non-trivial on $\Sigma$ or $\partial D_1$ and $\partial D_2$ are isotopic.

Next consider the case when the orientations at the intersections agree and $\Sigma \cap N$ is a twice punctured torus.  If the two loops in $\Sigma \cap \partial N$ are parallel then $\Sigma$ is the union of a twice punctured torus and an annulus.  This surface has genus two.  If both loops bound disks then the surface is a closed torus.  Because of the assumption that $\Sigma$ has genus at least three, we conclude that the induced automorphism is non-trivial.
\end{proof}

There exist strongly irreducible Heegaard splittings that have a pair of disks that intersect in two points, some of which will be constructed in Section~\ref{onesidesect}.

If $(\Sigma, H^-, H^+)$ is a weakly reducible Heegaard splitting with genus at least three, then we can always find a pair of disks that intersect in two points. Let $D_1 \subset H^-$, $D_2 \subset H^+$ be disks with disjoint boundaries.  If their boundaries are parallel then replace $D_1$ with a second disk in $H^-$ that is disjoint from the original $D_1$ (and therefore from $\partial D_2$) and not parallel to $D_1$.  

Choose an arc $\alpha$ with one endpoint in $\partial D_1$, the second endpoint in $\partial D_2$ and interior disjoint from both loops.  Dragging the boundary of $D_1$ along this arc produces a disk $D'_1$ that intersects $\partial D_2$ in exactly two points.  The disks $D'_1$, $D_2$ satisfy the conditions of Lemma~\ref{twopointslem}, so we conclude that the kernel of the homomorphism $i$ contains an infinite order, reducible element.

If $\Sigma$ is weakly reducible and genus two, we note that a weakly reducible, genus two Heegaard splitting is always reducible. (This is a straightforward exercise.) Let $D_1 \subset H^-$ be a separating disk whose boundary also bounds a disk in $H^+$.  Let $D_2 \subset H^+$ be a non-separating disk whose boundary is disjoint from $\partial D_1$.  As in the proof of Lemma~\ref{twopointslem}, we can isotope $\partial D_2$ to intersect $\partial D_1$ in two points.  Let $N$ be a regular neighborhood in $\Sigma$ of these two loops.  The boundary of $N$ consists of a trivial loop, two parallel loops and one non-trivial, separating loop.  Thus a composition of Dehn twists along these loops will not be trivial in $Mod(\Sigma)$, but trivial in $Mod(M)$. Thus we conclude the following:

\begin{Coro}
If $(\Sigma, H^-, H^+)$ is a weakly reducible Heegaard splitting then the kernel of $i$ contains a reducible element.
\end{Coro}

Note that in this construction, the orientations defined by the two intersections agree, so the automorphism of $\Sigma$ consists of Dehn twists along the boundary loops of a four punctured sphere.  However, one of these loops is trivial in $\Sigma$, so we are in fact twisting along the boundary loops of a pair of pants, specifically a regular neighborhood of $\partial D_1 \cup \alpha \cup \partial D_2$.

\section{Open book decompositions}
\label{opensect}

An open book decomposition for a 3-manifold $M$ is an ordered pair $(L, \phi)$ where $L \subset M$ is a link and $\phi : (M \setminus L) \rightarrow S^1$ is a fibration of the complement of $L$.  One can construct a Heegaard splitting for $M$ from the open book decomposition $(L,\phi)$ as follows:  

Parametrize $S^1$ as the interval $[0,1]$ with its endpoints identified.  The preimage $\phi^{-1}((0,\frac{1}{2}))$ is homeomorphic to the interior of a punctured surface cross an interval.  Its closure in $M$ is a handlebody $H^-$ which is the union of $L$ and $\phi^{-1}([0,\frac{1}{2}])$.  The closure of the set $\phi^{-1}((\frac{1}{2},1))$ is a second handlebody $H^+$ such that $H^- \cap H^+$ is a surface $\Sigma$ which is the union of $L$ and the surfaces $\phi^{-1}(\{0,\frac{1}{2}\})$.  The triple $(\Sigma, H^-, H^+)$ is the Heegaard splitting for $M$ induced from $(L, \phi)$.  

\begin{Lem}
If $(\Sigma, H^-, H^+)$ is induced by an open book decomposition $(L, \phi)$ of $M$ then the kernel of $i$ contains a reducible element.
\end{Lem}

\begin{proof}
Let $N$ be a closed regular neighborhood of the link $L$.  The complement of $N$ is homeomorphic to the result of gluing the surface product $S \times [0,1]$ by a map $\psi : (S \times \{0\}) \rightarrow (S \times \{1\})$.

There is a family $\{h_t\}$ of maps from the complement of $N$ to itself such that $h_0$ is the identity and $h_t$ sends each leaf $S \times \{x\}$ onto $S \times \{x + t\}$ where $x + t$ is taken modulo 1.  Each $h_t$ extends to an automorphism of $M$ by rotating the solid tori of $N$.

For each integer $n$, $h_n$ takes the surface $\Sigma$ to itself and each handlebody onto itself.  The induced automorphism on $\Sigma$ fixes the loops $L \subset \Sigma$ and restricts to $\psi^n$ on the complement of $L$.
\end{proof}

Note that all Heegaard splittings induced by open book decompositions have distance at most two: If $\Sigma$ is induced by an open book decomposition, with leaf $S$, let $\alpha$ be a properly embedded, essential arc in $S$ and $\beta$ a disjoint, essential, simple closed curve in $S$.

The Heegaard splitting is determined by the handlebodies $S \times [0,\frac{1}{2}]$ and $S \times [\frac{1}{2}, 1]$ in which $\alpha \times [0,\frac{1}{2}]$ and $\alpha \times [\frac{1}{2}, 1]$, respectively, are properly embedded, essential disks.  Both disks are disjoint from the loop $\beta \times \{\frac{1}{2}\}$, which is essential in the Heegaard surface.  Thus the splitting has distance at most two.  

Conversely, Masur and Schleimer have shown (Lemma 12.12 in~\cite{ms:holes}) that the set of arcs in $S$ defined by the boundary of a disk $D$ in $S \times [0,\frac{1}{2}]$  by taking the intersection of $\partial D$ with $S \times \{0\}$ and $S \times \{\frac{1}{2}\}$ form a set with diameter at most 6 in the arc complex. If there are disjoint disks in $S \times [0,\frac{1}{2}]$ and $S \times [\frac{1}{2}, 1]$ then their arcs of intersection form a set of diameter at most 13 in the arc complex such that the arcs in $S \times \{0\}$ are disjoint from the images of the arcs in $S \times \{1\}$ under the monodromy. Thus if we choose a monodromy map that moves every arc in the arc complex for $S$ a distance of at least 15, then the induced Heegaard splitting will be strongly irreducible. Such a monodromy can be found, for example, by taking a high power of a pseudo-Anosov map.

\section{One-sided Heegaard splittings}
\label{onesidesect}

A \textit{one-sided Heegaard surface} is an embedded, non-orientable surface $S \subset M$ such that the complement $S \setminus M$ is an open handlebody $H$. A regular neighborhood of $S$ will be a twisted interval bundle such that $S$ is the surface that consists of the midpoints of the intervals in the interval bundle structure.

The mapping class groups of certain one-sided Heegaard surfaces were examined by Birman and Rubinstein~\cite{rb:homeo}. We can construct a (two-sided) Heegaard splitting from a one-sided Heegaard surface as follows: Let $N$ be a regular neighborhood of $S$ and let $\alpha$ be an interval in the interval bundle structure for $N$.  The complement in $N$ of a regular neighborhood $O \subset N$ of $\alpha$ is a handlebody $H^-$ whose genus is equal to the crosscap number of $S$. If $H$ is the complement $M \setminus N$ then the union $H \cup O$ will be a second handlebody $H^+$ and the common boundary $H^- \cap H^+$ will be a Heegaard surface $\Sigma$.

Note that if we isotope $\Sigma$ into $N$, it will form a strongly irreducible Heegaard surface for this interval bundle. The manifold $M$ is the result of gluing a handlebody to $N$, and Li has shown~\cite{li:amalg} that for a sufficiently complicated gluing map, such a Heegaard surface $\Sigma$ will be strongly irreducible. Thus we can use this construction to produce strongly irreducible Heegaard splittings.

An automorphism for $\Sigma$ can be induced by a one-sided Heegaard splitting in two different ways: First, if $\psi$ is an automorphism of $S$, i.e. a map from $M$ to itself that takes $S$ onto itself, then we can isotope the map $\psi$ to take intervals to intervals within $N$.  Moreover, because $\alpha$ is a single interval, we can choose $\psi$ so that it takes $\alpha$ onto itself as well as taking $O$ onto itself.  Thus $\psi$ will determine an element of $Mod(M, \Sigma)$.  If this element does not reverse the direction of $\alpha$ then it will fix a meridian dual to the arc $\alpha$ and thus be reducible.

To find a one-sided splitting with such a symmetry, let $L \subset M$ be a knot such that $M \setminus L$ is fibered and the boundary of any page of the fibration wraps around $L$ exactly twice. Then the union of $L$ and any page of the fibration forms a one-sided Heegaard surface $S$ for $M$ and the open book structure determines an automorphism of $S$ similar to the automorphisms of two-sided Heegaard surfaces induced by open books.

Assume that the arc used to turn $S$ into a two-sided surface $\Sigma$ lies in a regular neighborhood of $L$. Then outside this regular neighborhood, $\Sigma$ is very similar to the Heegaard splitting induced by an open book decomposition.  Inside the solid torus neighborhood $\Sigma$ looks like a twice-punctured torus defining a torus twist, as described in Section~\ref{wrsect}. Applying Masur and Schleimer's Lemma~\cite{ms:holes}, as we did in the open book construction, one can show that if the monodromy map for the surface bundle moves every arc in the arc complex a distance of at least 15 then this Heegaard splitting will have distance exactly two.

Whether or not the mapping class group of the one-sided surface $S$ is non-trivial, there is a second type of automorphism of $\Sigma$ as follows: Consider the fundamental group of $S$ with the point of intersection $S \cap \alpha$ as its base point.  Every element of this group defines a path, which we can extend to an isotopy of $\alpha$ within $N$.  This extends to an isotopy of $\Sigma$, which determines an element of $Mod(M, \Sigma)$. If the path in $S$ is orientation preserving then the automorphism will fix the meridian of $\alpha$. Thus we have the following:

\begin{Lem}
There exist strongly irreducible Heegaard splittings induced by one-sided Heegaard splitting and every Heegaard splitting induced by a one-sided Heegaard surface admits a reducible automorphism.
\end{Lem}

Note that if one carries out this construction starting with a two-sided Heegaard splittings, the result will always be stabilized. (In fact, it will produce a stabilization of the initial Heegaard splitting.) However, there is still a way to generalize this construction. 

Let $G \subset M$ be an embedded handlebody or a union of disjoint, embedded handlebodies and let $S \subset M$ be a connected (one-sided or two-sided) surface properly embedded in the complement of $G$ such that the boundary of $S$ is disk busting in $G$.  In other words, assume that for every essential, properly embedded disk $D \subset G$, there is a loop in $\partial S$ that cannot be isotoped disjoint from $D$. Moreover, assume that the complement of $S \cup G$ in $M$ is one or two handlebodies, depending on whether $S$ is one-sided or two-sided. We will call the structure $(G, S)$ a \textit{branched Heegaard surface}.

We can form a Heegaard splitting $(\Sigma, H^-, H^+)$ by a construction similar to the one for one-sided surfaces: Let $H^-$ be a regular neighborhood of $G$ union the result of puncturing $S$ in a single point.  Let $N$ be a regular neighborhood of $G \cup S$. In the case when $S$ is non-orientable, the complement of $N$ is a single handlebody and Li's Theorem~\cite{li:amalg} implies that there are examples where the resulting Heegaard splitting is strongly irreducible. 

In the case when $S$ is orientable, we form $M$ by gluing two handlebodies to the boundary of $N$, so we cannot apply Li's Theorem. Instead, let $O \subset N$ be a regular neighborhood of the arc $\alpha$ dual to the puncture in $S$. Choose a component $F \subset \partial N$ and let $E \subset F$ be the disk $O \cap F$. By construction, $N \setminus O$ is a handlebody $H'$ and $F' = F \setminus O$ is a subsurface of $\partial H'$. This subsurface is not the horizontal boundary of an interval bundle in $H'$, so Theorem~12.1 in~\cite{ms:holes} implies that the projection of the disk set for $H'$ into the arc complex for $F'$ has diameter at most 60.

Any arc in $F'$ can be extended uniquely into $E$ to form a simple closed curve.  If two arcs are disjoint then their corresponding loops will have distance at most two. Thus the diameter of the set of loops in $F$ defined by disks in $H'$ is at most 120. If we glue a handlebody $H''$ to $F$ so that the disk set of $H''$ is distance at least three from the projection of the disk set of $H'$ in the curve complex for $F$ then every disk in $H''$ will intersect every disk in $H'$, so the induced Heegaard splitting of the resulting 3-manifold will be strongly irreducible. By Li's Theorem~\cite{li:amalg} we can then glue in the second handlebody so that the final Heegaard splitting $(\Sigma, H^-, H^+)$ is strongly irreducible.

We have shown that strongly irreducible Heegaard splittings can result from branched Heegaard splittings with orientable or non-orientable surfaces. As with Heegaard splittings induced by one-sided Heegaard surfaces, there are reducible elements of $Mod(M, \Sigma)$ that result from isotoping the arc $\alpha$ around a path in $S$. Thus we have the following.

\begin{Lem}
There exist strongly irreducible Heegaard splittings induced by branched Heegaard splittings and every Heegaard splitting induced by a branched Heegaard splitting admits a reducible automorphism.
\end{Lem}

\section{Classifying reducible automorphisms}
\label{refsect}

\begin{proof}[Proof of Theorem~\ref{redthm}]
Let $\Sigma$ be a strongly irreducible Heegaard surface for a hyperbolic 3-manifold $M$ and $\phi \in Mod(M, \Sigma)$ an infinite order, reducible automorphism with fixed set $L \subset \Sigma$. Choose a power of $\phi$ that is the identity on each component of $L$.

Because $M$ is irreducible, if $M \setminus L$ is reducible then there is a sphere $S \subset M \setminus L$ bounding a ball $B \subset M$ containing one or more components of $L$. Choose this sphere to be transverse to $\Sigma$ and to intersect $\Sigma$ in the smallest possible number of loops among all such spheres.   If $S$ is compressible in the complement of $\Sigma$ then the compression produces two spheres that intersect $\Sigma$ in fewer loops, at least one of which is essential in the complement of $L$.

Thus $S$ must be incompressible in the complement of $\Sigma$. By Scharlemann's characterization of how a strongly irreducible Heegaard splitting intersects a ball (Theorem 2.1 in~\cite{schar:ball}), we must have that $\Sigma \cap B$ is an unknotted punctured sphere.  In this case, every component $\ell$ of $L$ in $B$ bounds a disk $D$ in one of the handlebodies for the Heegaard splitting, say $H^-$.  By Oertel's classification of reducible automorphisms of handlebodies~\cite{oer:hbdy}, either $\ell$ cobounds an essential annulus in $H^+$ or $\ell$ is in the vertical boundary of an interval bundle structure within $H^+$. 

In the first case, boundary compressing the annulus produces a disk in $H^+$ whose boundary is disjoint from $\ell = \partial D$.  Thus the first case cannot occur when $\Sigma$ is strongly irreducible. 

In the second case, let $X \subset H^+$ be the interval bundle and let $S$ be union of the surface formed by the midpoints of the intervals and the disk $D \subset H^-$ bounded by $\ell$. The closure of the complement $H^+ \setminus X$ is a (possibly empty) handlebody $G$ such that $\partial S$ is contained in $\partial G$. Any disk in $G$ disjoint from $\partial S$ is properly embedded in $H^+$ and disjoint from $D$. Because $\Sigma$ is strongly irreducible, this implies that $\partial S$ must be disk busting in $G$. Moreover, the complement of $G \cup S$ is isotopic to the handlebody that results from compressing $H^-$ across $D$. If $G$ is empty then $S$ is a one-sided or two-sided Heegaard surface, though it must be one-sided since $\Sigma$ is strongly irreducible. If $G$ is not empty then $(G, S)$ is a branched Heegaard surface.

Otherwise, we will assume $L$ is irreducible. If the complement of $L$ is not hyperbolic then there is a collection of JSJ tori $\mathcal{T} \subset M \setminus L$.  (Note that since $M$ is hyperbolic, $M \setminus L$ cannot be a Seifert fibred space with two exceptional fibres). Because $M$ is hyperbolic and $L$ is irreducible, each torus $T \subset \mathcal{T}$ bounds a solid torus $C$ in $M$. Again, we will compress $T$ maximally in the complement of $\Sigma$. By Scharlemann's classification of how a strongly irreducible Heegaard splitting intersects a handlebody~\cite{schar:hbdy}, the intersection $\Sigma \cap C$ consists of incompressible annuli and possibly one piece that results from attaching a tube between two such annuli, or from an annulus to itself.  (These are the types of pieces that determine torus twists.) By Oertel's classification~\cite{oer:hbdy}, the loops of $L$ must cobound annuli with each other or with components of $L$ outside of $T$, though they are not parallel in $\Sigma$. Thus the loops must all be parallel in $C$ and thus parallel into $T$.

Let $M' \subset M$ be the component of $M \setminus T$ that is not a solid torus $C$.  By construction, $M'$ is atoroidal and thus hyperbolic.  Because $\phi$ fixes each component of $L$ and $\mathcal{T}$ consists of a union of annuli in the handlebodies $H^-$, $H^+$ with boundary parallel into $L$, the map $\phi$ can be isotoped to fix $\mathcal{T}$, and thus take $M'$ onto itself.  

If $\Sigma$ is incompressible in the complement of $L$ then let $R = \Sigma$. Otherwise, there is either a single component of $\Sigma \setminus L$ that has compressing disks disjoint from $L$ on both sides, or there are one or more components with compressing disks all on the same side.  In the first case, let $R$ be the union of the two surfaces that result from compressing the compressible component of $\Sigma \setminus L$ maximally in each direction. In the second case, let $R$ be the result of compressing $\Sigma$ maximally in the complement of $L$ (all in one direction.) We can conclude that $R$ is incompressible, using  Scharlemann's no nesting lemma for compressing disks for strongly irreducible Heegaard splittings \cite{schar:ball}.

The surface $R$ is canonical up to isotopy, i.e.\ its isotopy class does not depend on the choice of compressions, as long as they are maximal. Thus $\phi$ can be isotoped to take $R$ onto itself.

The hyperbolic submanifold $M' \subset M$ has a finite mapping class group, so some power $\phi^k$ is isotopy trivial on $M'$ but restricts to a non-trivial automorphism of $\Sigma \cap M'$, since $\phi$ has infinite order.   This implies that there is an ambient isotopy of $M'$ that takes $R$ onto itself.  If the induced automorphism of $R$ from this isotopy is non-trivial then the isotopy of $M'$ that takes $R$ onto itself defines a map of a non-trivial surface bundle over $R$ into $M'$. Because $R$ is incompressible, the inclusion map of its fundamental group is injective, so this implies that the fundamental group of $M'$ contains a semi-direct product of the fundamental group of $R$ with $\mathbf{Z}$, since no power of the image loop $\alpha$ of a base point in $R$ under the isotopy can be homotoped into $R$, as the isotopy induces a non-trivial automorphism of $R$. (It is simple to arrange that the isotopy moves a base point, chosen on $R$, around a loop $\alpha$ in $M'$).

By the proof of Theorem 11.1 in~\cite{hem:book}, this implies $M'$ is a surface bundle with fiber $R$. Notice that without loss of generality, a boundary loop $C^\pm$ of $R$ is moved by the isotopy around a mapped in torus $T^\pm$. This torus is $\pi_1$-injective and since $M'$ is hyperbolic, the torus is homotopic to a covering of a boundary component of $M'$.

If $M$ is not a surface bundle, then there is an infinite cyclic covering of $M$ given by using the cohomology class dual to $R$, which evaluates $\alpha$ to a non-zero integer. The lifts of $R$ to this covering space bound regions which are not products and so the inclusions of their fundamental groups are not conjugate subgroups. We conclude that this cannot occur if there is an ambient isotopy taking $R$ to itself.

Consequently, $R$ defines a Heegaard surface $\Sigma'$ for $M$, after we extend it into the solid tori by annuli, or a twice-punctured torus if the slope wraps around a solid torus twice. In the latter case, if there are any components of $L$ outside of this solid torus, the surface $\Sigma'$ will be weakly reducible. Thus $\Sigma'$ comes from one of the two constructions described in Section~\ref{opensect}. The surface $\Sigma$ is the result of attaching zero or more tubes to $\Sigma'$. By Proposition 22 in~\cite{me:stabs}, attaching tubes in this way always produces a stabilized (and thus weakly reducible) Heegaard splitting. Thus $\Sigma$ must be isotopic to $\Sigma'$ and is induced by an open book decomposition.

Otherwise, if $\phi$ does restrict to the identity map on $R$ then the isotopy defined by $\phi$ takes place entirely in the complement of a regular neighborhood of $R$, i.e. in a submanifold of $M$ bounded by a union of surfaces that result from compressing subsurfaces of $\Sigma$ in the complement of $L$. This corresponds to the final case in the theorem. 
\end{proof}

\section{Stabilized Splittings}
\label{stabsect}

In this final section, we construct a large collection of pseudo Anosov automorphisms of any stabilized Heegaard splitting of genus three or more.  The same methods can be used to find pseudo Anosov automorphisms of stabilized genus two splittings, but we leave the details of this construction to the reader.  The pseudo Anosov automorphisms are compositions of reducible automorphisms.  

It is an open question whether there are Heegaard splittings with only pseudo Anosov automorphisms.  Such a splitting would need to be strongly irreducible.  Pseudo Anosov automorphisms do occur in irreducible Heegaard splittings, but there are no general methods known for constructing such splittings.  The only known example was found by Long~\cite{long:twoh}.

\begin{Thm}
If $(\Sigma, H^-, H^+)$ is a stabilized Heegaard splitting of genus at least three then $Mod(M, \Sigma)$ contains a pseudo Anosov element in the kernel of $i$.
\end{Thm}

\begin{proof}
Let $(\Sigma, H^-, H^+)$ be a stabilized Heegaard splitting of a manifold $M$.  Then there is an embedded sphere $S \subset M$ bounding a ball $B$ such that $S \cap \Sigma$ is a single loop and $\Sigma \cap B$ is a punctured torus. The sphere $S$ is called a \textit{reducing sphere}. If we replace the subsurface $\Sigma \cap B$ with one of the disks in the complement $S \setminus \Sigma$, the result is a Heegaard surface $R$ of a Heegaard splitting $(R, G^-, G^+)$.  

We can reconstruct $\Sigma$ from $R$ by choosing a ball $B \subset M$ such that the intersection $B \cap \Sigma$ is a disk $D \subset R$, then replacing $D$ with a standardly embedded punctured torus in $B$.  The sphere $\partial B$ becomes a reducing sphere for the new surface (which is isotopic to $\Sigma$.)  

Assume there is a properly embedded, essential disk $D_1 \subset H^-$ such that the intersection $D \cap \partial D_1$ consists of a single arc.  Let $\alpha \subset D_1$ be a properly embedded arc whose endpoints are in $D$.  Stabilizing $R$ at $D$ is equivalent to attaching a tube to $R$ along $\alpha$.  The reducing sphere $\partial B$ is isotopic to the union of the boundary of a regular neighborhood of $\alpha$ and one of the disk components of $D_1 \setminus \alpha$.  There is a second reducing sphere defined by the union of a regular neighborhood of $\alpha$ and the second component of $D_1 \setminus \alpha$.  The loops of intersection between $\Sigma$ and these spheres are shown in Figure~\ref{diskstabfig}.  As shown on the right of the figure, the loop determined by the second reducing sphere intersects $\Sigma \cap B$ in two arcs and intersects $R$ in a pair of arcs parallel to $\partial D_1 \setminus D$.
\begin{figure}[htb]
  \begin{center}
  \includegraphics[width=3in]{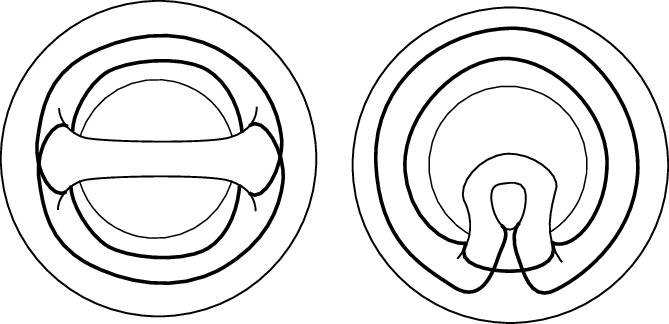}
  \caption{After adding a stabilization along an arc in a disk $D$ in $H^-$, there is a second reducing sphere defined by $D$.}
  \label{diskstabfig}
  \end{center}
\end{figure}

Because $\Sigma$ has genus at least three, $R$ has genus at least two.  Hempel~\cite{hempel} showed that the diameter of the set of curves in the curve complex bounding loops in a handlebody of genus at least two is infinite.  Thus there are disks $D'_1 \subset G^-$ and $D'_2 \subset G^+$ whose boundaries are distance at least three apart.  This implies that the complement in $R$ of the boundaries of these disks is a collection of disks.
\begin{figure}[htb]
 \begin{center}
  \includegraphics[width=3.5in]{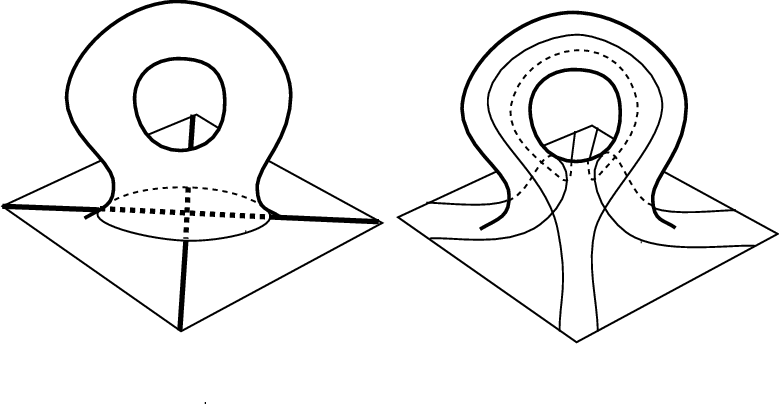}
  \caption{By placing the stabilization at the intersection of a pair of disks in opposite handlebodies, we find two reducing spheres that cut the surface into disks.}
  \label{stabfig}
  \end{center}
\end{figure}

Let $x$ be a point in the intersection $\partial D'_1 \cap \partial D'_2$ and let $D \subset R$ be a small disk containing $x$ and such that $D$ intersects the boundary of each disk $D'_1$, $D'_2$ in a single arc.  Construct a surface isotopic to $\Sigma$ by stabilising $\Sigma$ at the disk $D$.

As described above, there is a reducing sphere $S_1$ for $\Sigma$ such that $S_1 \cap R$ consists of two arcs parallel to $\partial D_1$ and a second reducing sphere $S_2$ that intersects $R$ in a pair of arcs parallel to $\partial D_2$, as shown in Figure~\ref{stabfig}.  The complement in $R$ of $S_1 \cup S_2$ is a collection of disks.  The complement in $\Sigma \cap B$ of $S_1 \cup S_2$ consists of three disks, as in Figure~\ref{stabfig}.

There is an automorphism of $(\Sigma, H^-, H^+)$ that restricts to a Dehn twist along $S_1 \cap \Sigma$ and a second automorphism that restricts to a Dehn twist along $S_2 \cap \Sigma$.  Because the complement in $\Sigma$ of the loops $S_1 \cap \Sigma$ and $S_2 \cap \Sigma$ is a collection of disks, a Theorem of Penner~\cite{pen:pseudo} implies that a composition of powers of these two Dehn twists is pseudo Anosov.
\end{proof}

\bibliographystyle{amsplain}
\bibliography{finitemcg}

\providecommand{\bysame}{\leavevmode\hbox to3em{\hrulefill}\thinspace}
\providecommand{\MR}{\relax\ifhmode\unskip\space\fi MR }
\providecommand{\MRhref}[2]{%
  \href{http://www.ams.org/mathscinet-getitem?mr=#1}{#2}
}
\providecommand{\href}[2]{#2}
\begin{thebibliography}{10}

\bibitem{akbas}
Erol Akbas, \emph{A presentation for the automorphisms of the 3-sphere that
  preserve a genus two {H}eegaard splitting}, Pacific J. Math. \textbf{236}
  (2008), no.~2, 201--222. \MR{2407105 (2009d:57029)}

\bibitem{bach}
David Bachman, \emph{Connected sums of unstabilized {H}eegaard splittings are
  unstabilized}, Geom. Topol. \textbf{12} (2008), no.~4, 2327--2378.
  \MR{2443968 (2009h:57035)}

\bibitem{bdt:seif}
David Bachman and Ryan Derby-Talbot, \emph{Non-isotopic {H}eegaard splittings
  of {S}eifert fibered spaces}, Algebr. Geom. Topol. \textbf{6} (2006),
  351--372 (electronic), With an appendix by R. Weidmann. \MR{2220681
  (2007a:57012)}

\bibitem{rathcam}
Marion~Moore Campisi and Matt Rathbun, \emph{High distance knots in closed
  3-manifolds}, preprint (2009), arXiv:0911.3675.

\bibitem{oc:auts}
Leonardo~N. Carvalho and Ulrich Oertel, \emph{A classification of automorphisms
  of compact 3-manifolds}, preprint (2005).

\bibitem{cl:alg}
A.~J. Casson and D.~D. Long, \emph{Algorithmic compression of surface
  automorphisms}, Invent. Math. \textbf{81} (1985), no.~2, 295--303. \MR{799268
  (86m:57012)}

\bibitem{cho}
Sangbum Cho, \emph{Homeomorphisms of the 3-sphere that preserve a {H}eegaard
  splitting of genus two}, Proc. Amer. Math. Soc. \textbf{136} (2008), no.~3,
  1113--1123 (electronic). \MR{2361888 (2009c:57029)}

\bibitem{mccho}
Sangbum Cho and Darryl McCullough, \emph{The tree of knot tunnels}, Geom.
  Topol. \textbf{13} (2009), no.~2, 769--815. \MR{2469530 (2010j:57005)}

\bibitem{dl:ricci}
Jonathan Dinkelbach and Bernhard Leeb, \emph{Equivariant {R}icci flow with
  surgery and applications to finite group actions on geometric 3-manifolds},
  Geom. Topol. \textbf{13} (2009), no.~2, 1129--1173. \MR{2491658
  (2011b:53158)}

\bibitem{fu:invol}
David Futer, \emph{Involutions of knots that fix unknotting tunnels}, J. Knot
  Theory Ramifications \textbf{16} (2007), no.~6, 741--748. \MR{2341313
  (2008k:57009)}

\bibitem{goer:s3}
L.~Goeritz, \emph{{Die Abbildungen der Brezelfl\"{a}che und der Volbrezel vom
  Gesschlet 2}}, Abh. Math. Sem. Univ. Hamburg \textbf{9} (1933), 244--259.

\bibitem{htt}
Joel Hass, Abigail Thompson, and William Thurston, \emph{Stabilization of
  {H}eegaard splittings}, Geom. Topol. \textbf{13} (2009), no.~4, 2029--2050.
  \MR{2507114 (2010k:57044)}

\bibitem{hem:book}
John Hempel, \emph{{$3$}-{M}anifolds}, Princeton University Press, Princeton,
  N. J., 1976, Ann. of Math. Studies, No. 86. \MR{0415619 (54 \#3702)}

\bibitem{hempel}
\bysame, \emph{3-manifolds as viewed from the curve complex}, Topology
  \textbf{40} (2001), no.~3, 631--657. \MR{1838999 (2002f:57044)}

\bibitem{me:t3}
Jesse Johnson, \emph{{Automorphisms of the three-torus preserving a genus three
  Heegaard splitting}}, preprint (2007), arXiv:0708.2683.

\bibitem{me:stabs}
\bysame, \emph{Bounding the stable genera of {H}eegaard splittings from below},
  J. Topol. \textbf{3} (2010), no.~3, 668--690. \MR{2684516 (2011g:57024)}

\bibitem{me:medium}
\bysame, \emph{Mapping class groups of medium distance {H}eegaard splittings},
  Proc. Amer. Math. Soc. \textbf{138} (2010), no.~12, 4529--4535. \MR{2680077}

\bibitem{JMc}
Jesse Johnson and Darryl McCullough, \emph{{The space of Heegaard splittings}},
  preprint (2010), arXiv:1011.0702.

\bibitem{ker:niels}
Steven~P. Kerckhoff, \emph{The {N}ielsen realization problem}, Ann. of Math.
  (2) \textbf{117} (1983), no.~2, 235--265. \MR{690845 (85e:32029)}

\bibitem{li:wald}
Tao Li, \emph{Heegaard surfaces and measured laminations. {I}. {T}he
  {W}aldhausen conjecture}, Invent. Math. \textbf{167} (2007), no.~1, 135--177.
  \MR{2264807 (2008h:57033)}

\bibitem{li:amalg}
\bysame, \emph{Heegaard surfaces and the distance of amalgamation}, Geom.
  Topol. \textbf{14} (2010), no.~4, 1871--1919. \MR{2680206}

\bibitem{long:twoh}
D.~D. Long, \emph{On pseudo-{A}nosov maps which extend over two handlebodies},
  Proc. Edinburgh Math. Soc. (2) \textbf{33} (1990), no.~2, 181--190.
  \MR{1057747 (91j:57016)}

\bibitem{lm:dblrect}
Martin Lustig and Yoav Moriah, \emph{A finiteness result for {H}eegaard
  splittings}, Topology \textbf{43} (2004), no.~5, 1165--1182. \MR{2079999
  (2005d:57031)}

\bibitem{ms:holes}
Howard Masur and Saul Schleimer, \emph{The geometry of the disk complex},
  preprint (2010), arXiv:1010.3174.

\bibitem{meeksyau}
William~H. Meeks, III and Shing~Tung Yau, \emph{The equivariant {D}ehn's lemma
  and loop theorem}, Comment. Math. Helv. \textbf{56} (1981), no.~2, 225--239.
  \MR{630952 (83b:57006)}

\bibitem{ms:tunnels}
Kanji Morimoto and Makoto Sakuma, \emph{On unknotting tunnels for knots}, Math.
  Ann. \textbf{289} (1991), no.~1, 143--167. \MR{1087243 (92e:57015)}

\bibitem{nam:mcg}
Hossein Namazi, \emph{Big {H}eegaard distance implies finite mapping class
  group}, Topology Appl. \textbf{154} (2007), no.~16, 2939--2949. \MR{2355879
  (2008j:57025)}

\bibitem{oer:hbdy}
Ulrich Oertel, \emph{Automorphisms of three-dimensional handlebodies}, Topology
  \textbf{41} (2002), no.~2, 363--410. \MR{1876895 (2002j:57039)}

\bibitem{pen:pseudo}
Robert~C. Penner, \emph{A construction of pseudo-{A}nosov homeomorphisms},
  Trans. Amer. Math. Soc. \textbf{310} (1988), no.~1, 179--197. \MR{930079
  (89k:57026)}

\bibitem{rb:homeo}
J.~H. Rubinstein and J.~S. Birman, \emph{One-sided {H}eegaard splittings and
  homeotopy groups of some {$3$}-manifolds}, Proc. London Math. Soc. (3)
  \textbf{49} (1984), no.~3, 517--536. \MR{759302 (85m:57009)}

\bibitem{schar:ball}
Martin Scharlemann, \emph{Local detection of strongly irreducible {H}eegaard
  splittings}, Topology Appl. \textbf{90} (1998), no.~1-3, 135--147.
  \MR{1648310 (99h:57040)}

\bibitem{schar:hbdy}
\bysame, \emph{Local detection of strongly irreducible {H}eegaard splittings},
  Topology Appl. \textbf{90} (1998), no.~1-3, 135--147. \MR{1648310
  (99h:57040)}

\bibitem{schar:mcg}
\bysame, \emph{Automorphisms of the 3-sphere that preserve a genus two
  {H}eegaard splitting}, Bol. Soc. Mat. Mexicana (3) \textbf{10} (2004),
  no.~Special Issue, 503--514. \MR{2199366 (2007c:57020)}

\bibitem{zim:Hur}
B.~Zimmermann, \emph{Hurwitz groups and finite group actions on hyperbolic
  {$3$}-manifolds}, J. London Math. Soc. (2) \textbf{52} (1995), no.~1,
  199--208. \MR{1345726 (96k:57011)}

\end{thebibliography}

\end{document}